\documentclass[11pt,twoside]{article}
\usepackage{amsmath,amsbsy,amsfonts}
\usepackage{color}
\newtheorem{theorem}{Theorem}[section]
\newtheorem{lemma}[theorem]{Lemma}
\newtheorem{proposition}[theorem]{Proposition}

\newtheorem{remark}[theorem]{Remark}

\newcommand{\qed}{\hfill\rule{2mm}{2mm}}

\title{A sign-changing solution for the Schr\"odinger-Poisson equation\footnote{Partially supported by INCT-MAT, Casadinho/PROCAD 552464/2011-2.}}

\author{Claudianor O. Alves\footnote{C.O. Alves was partially supported by CNPq/Brazil  301807/2013-2},\  Marco A. S. Souto\footnote{M.A.S. Souto was
supported by CNPq/Brazil 304652/2011-3.}\\
\noindent Unidade Acad\^emica de Matem\'atica e Estat\'istica\\
\noindent Universidade Federal de Campina Grande \\
\noindent 58109-970, Campina Grande - PB, Brazil.\\
\noindent e-mail: {\tt{coalves@dme.ufcg.edu.br, marco@dme.ufcg.edu.br}}\\
\mbox{} and
\mbox{}\\
\noindent S\'{e}rgio H. M. Soares\footnote{Corresponding author.
 E-mail:~\textsf{monari@icmc.usp.br}, Phone: +55\,16\,3373\,9660, Fax: +55\,16\,3373\,9650.}\\
\noindent Departamento de Matem\'atica \\
\noindent Instituto de Ci\^{e}ncias Matem\'{a}ticas e de
 Computa\c{c}\~{a}o\\
\noindent Universidade de S\~ao Paulo \\
\noindent 13560-970, S\~ao Carlos - SP, Brazil. \\
\noindent e-mail: {\tt{monari@icmc.usp.br}} }
\date{}

\begin{document}
\maketitle



\begin{abstract}

We find a sign-changing solution for a class of
Schr\"odinger-Poisson system in $\mathbb{R}^3$  as an existence result by minimization in a closed subset containing all the sign-changing solutions of the equation. The proof is based on variational methods in association with the deformation lemma and Miranda's theorem.
\end{abstract}

\section{Introduction}

This paper is concerned with establishing sign-changing for a class of nonlinear Schr\"{o}dinger-Poisson equations, which includes the typical and relevant model case
$$
\left\{ \begin{array}{ll}
  - \Delta u +V(x)u+
 \phi u = |u|^{p-2}u\  \mbox{ in $\mathbb R^3$},\\
-\Delta \phi=u^2\  \mbox{ in $\mathbb R^3$}.
\end{array}
\right.\eqno{(NLSP)}
$$
The Schr\"{o}dinger-Poisson equation has a great importance in the study of stationary solutions $\psi(x,t) = e^{-it}u(x)$ of  the time-dependent Schr\"odinger-Poisson equation, which describes quantum (non-relativistic) particles interacting with the eletromagnetic field generated by the motion (see \cite{AM,bf,BM,Mauser,Ruiz,Sanchez} for more details).

Many recent studies of (NLSP) have focused on existence and
nonexistence of solutions, multiplicity of solutions, ground states, radially
and non-radial solutions, semiclassical limit and concentrations of solutions
(see \cite{JMAA2010,nodal1,AM, ap, Cerami, C, Daprile1, Daprile2, Davenia, G, K1, Ianni, Ianni2, Kim, LG, Ruiz,Ruiz-Sic} and the references given there). In \cite{C}, Coclite proved the existence of a nontrivial radial solution of
(NLSP)  when $4 < p < 6$ and $V$ is a positive constant. The same result was
established in \cite{Daprile1} for $4 \leq p < 6$. In \cite{Daprile2}, by using a Pohozaev-type identity,
D'Aprile and Mugnai proved that (NLSP) has no nontrivial solution when $p \leq 2$ or
$p \geq 6$. This result was completed in \cite{Ruiz}, where Ruiz showed that if $p \leq 3$, the problem (NLSP) does not admit any nontrivial solution, and if  $3 < p < 6$, there exists a nontrivial radial solution of (NLSP). In \cite{ap}, Azzollini and Pomponio proved the
 existence of ground state solutions of (NLSP) when $3 < p < 6$ and $V$ is
a positive constant. The case non-constant potential was also treated in \cite{ap} for
$4 < p < 6$ and $V$ is possibly unbounded below.  The Schr\"odinger-Poisson equation has also been considered in a bounded domain in the papers of Siciliano \cite{G}, Ruiz and Siciliano \cite{Ruiz-Sic} and Pisani and Siciliano \cite{LG}. All these papers are about  positive solutions to (NLSP).  There are few results about sign-changing solutions to (NLSP).  The best references here are \cite{nodal1}, \cite{Ianni2}, and \cite{Kim}. Ianni \cite{Ianni2} employed a dynamical approach (not variational) in order to show the existence of radial solutions to (NLSP) {for $V$ constant} and $p\in [4, 6)$  having  a prescribed number of nodal domains. To obtain this result, she first studied the existence of sign changing radial solutions for the corresponding (NLSP) in balls of $\mathbb{R}^3$ with Dirichlet boundary conditions. Kim and Seok \cite{Kim} obtained results similar to \cite{Ianni2}  for $p \in (4,6)$ by using an extension of the Nehari variational method \cite{nehari, terracininodea}.
Alves and Souto \cite{nodal1} considered the problem (NLSP)  in a bounded domain and $V \equiv 0$ and proved the existence of least energy  sign-changing solutions for (NLSP)  changing sign exactly once in $\mathbb{R}^3$. The proof is based on variational methods. More precisely, it was proved that the associated energy functional assumes a minimum value on the nodal set.

Motivated by the results just described, we are interested in finding sign-changing
solution for (NLSP) in $\mathbb{R}^3$, {where potential $V$  is not necessarily a radially symmetric function}. The result will be state for a class of more general problem
$$
\left\{ \begin{array}{ll}
  - \Delta u +V(x)u+
 \phi u = f(u)  \mbox{ in $\mathbb R^3$},\\
-\Delta \phi=u^2  \mbox{ in $\mathbb R^3$},
\end{array}
\right.\eqno{(SP)}
$$
where $f$ belongs to $C^1(\mathbb R,\mathbb R)$ and satisfies
\begin{enumerate}
\item[$(f_1)$\ ]
$\displaystyle\lim_{s\rightarrow 0 } \frac{f(s)}{s} = 0$;
\item[$(f_2)$\ ] $\displaystyle\lim_{|s|\rightarrow +\infty } \frac{f(s)}{|s|^5} = 0$;
\item[$(f_3)$\ ] $\displaystyle\lim_{|s|\rightarrow +\infty} \frac{F(s)}{s^4}
=+\infty$;
\item[$(f_4)$\ ] $\displaystyle \frac{f(s)}{s^3} $ is increasing in $|s|>0$.
\end{enumerate}
\begin{remark}\label{rem} We observe that $(f_2)$ is weaker than the usual subcritical condition.  The conditions  $(f_1)$ and $(f_4)$   imply that $H(s)=sf(s)-4F(s)$ is a non-negative function, increasing  in $|s|$ with
$$
sH'(s)=s^{2}f'(s)-3f(s)s>0 \,\,\, \mbox{for any} \,\, |s|>0.
$$
Note that $(f_4)$ is a weaker version of the famous Ambrosetti-Rabinowitz condition for this class of equation. We need to suppose this condition because we will apply variational methods in our argument, and as we shall see, the associated energy functional has a local term of the fourth order.
\end{remark}

Here $V:\mathbb R^3 \to \mathbb R$  is locally H\"older continuous  and satisfies the following assumptions:
\begin{enumerate}
\item[$(V_1)$\ ]  There exists $\alpha >0$ such that  $V(x) \geq \alpha >0$, for all $x \in \mathbb R^3$;
\item[$(V_2)$\ ] $V_\infty=\max \{V(x):x\in\mathbb R^3\}$ and
\[
\lim_{|x|\rightarrow +\infty} V(x)=V_\infty;
\]
\item[$(V_3)$\ ]  There exist $R_o >0$  and    $\rho : (R_o, \infty) \to (0, \infty)$ a  non-decreasing
function such that
\[
\lim_{r\to \infty} \rho(r)e^{\delta r}=\infty
\]
for all $\delta >0$, and
$$ V(x)\leq V_\infty-\rho(|x|),  \mbox{
for all $|x|\geq R_o$}.
$$
\end{enumerate}

Our main result is the following:
\begin{theorem}{\label{fth}}
Suppose that   $f$  satisfies $(f_1)- (f_4)$ and $V$ satisfies $(V_1)-(V_3)$. Then problem ($SP$)
possesses a least energy sign-changing solution, which changes sign exactly once in $\mathbb{R}^3$.
\end{theorem}

Our result can be seen as a similar version for $\mathbb{R}^3$ of the result due to Alves and Souto \cite{nodal1}, which will be useful in our argument. We observe that Theorem \ref{fth} establishes the existence of a least energy sign-changing solution when  $f(s)= |s|^{p-2}s$, for  $p \in (4, 6)$.

As observed in \cite{BWW,BW}, the general procedures to find sign-changing solutions of an equation with a nonlinear term stumble in the fact that the nodal set is not a submanifold of $H^1$ because the map $u \mapsto u^\pm$ lacks differentiability; thus it is not evident that  a minimizer of the associated energy functional  on the nodal set  is a solution of the equation. Furthermore, there is a worsening in the case considered here: Since the associated energy functional has a nonlocal term, it follows that even if $u$ is a sign-changing solution of the problem, the functions $u^\pm$ do not belong both to the Nehari manifold, and so some arguments used to prove the existence of nodal solutions for semilinear local problems can not be used in our arguments.


Our approach is based on some arguments presented in \cite{nodal1, BWW} in association with the deformation lemma and Miranda's theorem. The contribution of our work are twofold: On one hand, it is applying the  construction of \cite{nodal1} in an unbounded domain like $\mathbb{R}^3$ and consequently dealing with the difficulties it brings; on the other hand, facing the subtle peculiarities of a nonlocal term.
We start by establishing some estimates involving functions that change sign.  We find a sign-changing solution as an existence result by minimization in a closed subset containing all the sign-changing solutions of the equation. At first, this may resemble the ideas found in \cite{BWW,BW}. However, we need to choose a suitable  minimizing sequence to the nodal level. This choice involves the corresponding equation in bounded domains (balls) and
the problem is then proving that the minimum of the energy on the corresponding closed subset containing all the sign-changing solutions of the equation in bounded domains is achieved by some function in the subset. In order to overcome the possible lack of regularity  of this subset, it is crucial to apply a deformation lemma and a fine use of Miranda's theorem \cite{Miranda}.

\section{The variational framework}

In this section, we present the variational framework to deal with problem $(SP)$. The key observation is that equation $(SP)$ can be transformed into a Schr\"{o}dinger equation  with a nonlocal term (see, for instance, \cite{ap,G,Ruiz}). This permits the use of variational methods. Effectively, by the Lax-Milgram Theorem,
given $u\in H^1(\mathbb R^3)$, there exists a unique $\phi=\phi_u
\in D^{1,2}(\mathbb R^3)$ such that
$$
-\Delta \phi=u^2.
$$
By using standard arguments, we have that $\phi_u$ verifies the following properties (for a proof see
\cite{Daprile1,Ruiz}):
\begin{lemma}{\label{lm1}} For any $u\in H^1(\mathbb R^3)$, we have
\begin{itemize}
\item[$i)$]  there exists $C>0$ such that
$||\phi_u||_{D^{1,2}}\leq C||u||_{H^1}^2$ and
$$
\int_{\mathbb R^3}|\nabla \phi_u|^2dx=\int_{\mathbb R^3}\phi_u
u^2dx\leq C ||u||_{H^1}^4 \quad \forall\, u\in H^1(\mathbb R^3);
$$
where $||u||_{H^1}^2=\int_{\mathbb R^3}(|\nabla u|^2+u^2)dx$ and
$||w||_{D^{1,2}}^2=\int_{\mathbb R^3}|\nabla w|^2dx$

\item[$ii)$] $\phi_u\geq 0$ $\forall u\in H^1(\mathbb R^3)$;

\item[$iii)$] $\phi_{tu}=t^2\phi_u$, $\forall t>0$ and $u\in H^1(\mathbb R^3)$;

\item[$iv)$] if $a \in \mathbb R^3$ and $u_a(x)=u(x-a)$ then
$\phi_{u_a}(x)=\phi_u(x-a)$ and
\[
\int_{\mathbb R^3}\phi_{ u_a}  u_a ^2dx=\int_{\mathbb R^3}\phi_{ u}
u ^2dx;
\]
\item[$v)$] if $u_n \rightharpoonup u$ in $H^1(\mathbb R^3)$, then
$\phi_{u_n} \rightharpoonup \phi_u$ in $H^1(\mathbb R^3)$ and
$$
\liminf_{n \rightarrow +\infty} \int_{\mathbb R^3}\phi_{ u_n}  u_n
^2dx\leq\int_{\mathbb R^3}\phi_u u^2dx.
$$
\end{itemize}
\end{lemma}
Therefore, $(u,\phi) \in H^1(\mathbb R^3)\times D^{1,2}(\mathbb
R^3)$ is a solution of $(SP)$ if, and only if, $\phi = \phi_u$ and
$u \in H^1(\mathbb R^3)$ is a weak solution of the nonlocal problem
$$
\left\{ \begin{array}{l}
  - \Delta u +V(x)u+ \phi_u u = f(u)  \mbox{ in }\mathbb R^3,\\
u\in H^1( \mathbb R^3) .
\end{array}
\right.\eqno{(P)}
$$
Combining ($f_1$)-($f_2$) with  Lemma \ref{lm1}, the  functional $J: H^1(\mathbb
R^3)\rightarrow \mathbb R$  given by
$$
J(u)=\frac 12 ||u||^2 +\frac 14\int_{\mathbb R^3}\phi_u u^2dx-
\int_{\mathbb R^3}F(u)dx,
$$
where
$$
||u||^2=\int_{\mathbb R^3}(|\nabla u|^2+V(x)u^2)dx, \quad F(s)=\int_0^s
f(t)dt,
$$
belongs to $C^1(H^{1}(\mathbb R^3),\mathbb{R})$ and
\[
J'(u)v =  \int_{\mathbb R^3}(\nabla u \nabla v +V(x)uv)dx +
\int_{\mathbb R^3}\phi_u uvdx - \int_{\mathbb R^3}f(u)vdx \,\,\,\,
\forall u,v \in H^1(\mathbb R^3).
\]
Hence, critical points of $J$ are the weak solutions for nonlocal problem $(P)$.

In the proof of Theorem \ref{fth}, we prove that functional $J$ assumes
a minimum value on the nodal set
$$
\mathcal M=\{u\in \mathcal N: J'(u)u^+=J'(u)u^-=0 \mbox { and }
u^{\pm}\neq 0 \}
$$
where $u^{+}=\max\{u(x),0\}$, $u^{-}(x)=\min\{u(x),0\}$ and
$$
\mathcal{N} = \{ u\in H^1(\mathbb R^3)\setminus\{0\}\, :\,
J'(u)u=0\}
$$
is the Nehari manifold associated with $J$. More precisely, we prove that there is $w \in \mathcal{M}$ such that
$$
J(w)=c_0\doteq\inf_{u \in \mathcal{M}}J(u).
$$
Furthermore, we prove that $w$ is a critical point of $J$, and so, $w$ is
a least energy nodal solution for $(P)$ with exactly two nodal
domains.

Since  $J$ has the nonlocal term $\int_{\mathbb R^3}\phi_u u^2dx$,
if $u$ is a sign-changing solution for $(P)$, we have that
$$
J'(u^+)u^+ =-\int_{\mathbb R^3}\phi_{u^{-}}(u^{+})^{2} <0\,\,\,
\mbox{and} \,\,\, J'(u^-)u^- =-\int_{\mathbb
R^3}\phi_{u^{+}}(u^{-})^{2} <0.
$$
Consequently, even though $u$ was a sign-changing solution for $(P)$, the functions $u^\pm$ do not belong both to $\mathcal{N}$. Hence, some arguments used to prove the existence of sign-changing  solutions for problem like
$$
\left\{ \begin{array}{l}
  - \Delta u +u= f(u)  \mbox{ in }\mathbb R^3,\\
u\in H^1( \mathbb R^3)
\end{array}
\right.\eqno{(P_1)}
$$
can not be applied; thus a careful analysis  is necessary in a lot of estimates.


\section {Technical lemmas }

{Consider the Sobolev space $H^1(\mathbb R^3)$ endowed with the norm
$$
||u||_*^2=\int_{\mathbb R^3}(|\nabla u|^2+V_\infty u^2)dx.
$$
Let $J_\infty: H^1(\mathbb R^3) \to \mathbb R$ be the functional given by
$$
J_\infty(u)=\frac 12 ||u||_*^2 +\frac
14\int_{\mathbb R^3}\phi_u u^2dx- \int_{\mathbb R^3}F(u)dx,
$$}
and consider
$$
\mathcal{N}_\infty = \{ u\in H^1(\mathbb R^3)\setminus\{0\}\, :\,
J_\infty'(u)u=0\},\quad  c_\infty= \inf_{\mathcal{N}_\infty} J_\infty.
$$

The following lemma establishes that {$c_0$} is a positive level. Similar result holds for $c_\infty$.

\begin{lemma}\label{lema2}
There exists $\rho>0$ such that
\begin{enumerate}
\item [(i)] $J(u)\geq ||u||^2/4$ and $||u||\geq \rho, \forall\, u\in \mathcal N$;
\item [(ii)] $||w^\pm||\geq\rho, \,\, \forall\, w \in \mathcal M$.
\end{enumerate}
\end{lemma}
\noindent {\bf{Proof:}} From Remark \ref{rem}, for every $u\in \mathcal N$,
\[
4J(u) = 4J(u) - J'(u)u = \|u\|^2 + \int_{\mathbb{R}^3}(f(u)u - 4F(u))dx \geq \|u\|^2,
\]
and (i) follows.  Taking $\alpha >0$, given by $(V_1)$, we set $\epsilon \in (0, \alpha)$. Since $f$ satisfies $(f_1)-(f_2)$,  there exists $C=C(\epsilon) >0$ such that
\begin{equation}\label{estimativa}
f(s)s \leq \epsilon s^2 + Cs^6, \forall\, s \in \mathbb{R}.
\end{equation}
For every  $w \in \mathcal M$, we have $J'(w^\pm)w^\pm <0$, which gives
\[
\|w^\pm\|^2 \leq \|w^\pm\|^2 + \int_{\mathbb{R}^3}\phi_{w^\pm}(w^\pm)^2 dx < \int_{\mathbb{R}^3} f(w^\pm)w^\pm dx.
\]
From (\ref{estimativa}), we have
\begin{eqnarray*}
\|w^\pm\|^2 &\leq& \epsilon \int_{\mathbb{R}^3}(w^\pm)^2dx + C\int_{\mathbb{R}^3}(w^\pm)^6dx \\
        &\leq& \frac{\epsilon}{\alpha}\int_{\mathbb{R}^3}V(x)(w^\pm)^2dx + C\int_{\mathbb{R}^3}(w^\pm)^6dx \\
         &\leq& \frac{\epsilon}{\alpha}\|w^\pm\|^2 + C\|w^\pm\|^6,
\end{eqnarray*}
and (ii) is proved. \qed

%
\medskip
The following lemma is a consequence of Miranda's theorem. A proof can be found in \cite{nodal1}.

\begin{lemma}\label{lema3}
Let  $v\in H^1(\mathbb R^3)$ satisfy $v^{\pm}\neq 0$. Then, there are
$t,s>0$ such that $J'(tv^++sv^-)v^+=0$ and $J'(tv^++sv^-)v^-=0$.
Moreover, if
$
J'(v)(v^{\pm}) \leq 0,
$
we have $s,t\leq 1$.
\end{lemma}

\section{The choice of the minimizing sequence}

Given $R>0$, let $B_R$ be the ball of radius $R$ centered at $0$. Consider the problem
$$
\left\{\begin{array}{l}
- \Delta u +V(x)u+\phi u= f(u),  \mbox{ in }  B_R, \\
-\Delta \phi=\tilde{u}^2,\mbox{ in }  \mathbb R^3,\\
 \phi\in D^{1,2}(\mathbb R^3), u \in H_0^1(B_R),
\end{array}\right. \leqno (AP_R)
$$
{where
$$
\tilde{u}(x)=
\left\{
\begin{array}{cl}
u(x) & \mbox{if $x\in B_R$}, \\
0  & \mbox{if $x\in \mathbb{R}^3\setminus B_R$}.
\end{array}
\right.
$$
{By Proposition \ref{appendix}, in Appendix, for any $R>0$, there exists a sign-changing solution $u=u_R$ of  $(AP_R)$  such that
\[
c_R= \inf_{u\in \mathcal{M}_R}J_R(u)=  J_R(u_R),
\]
where $J_R:H^{1}_{0}(B_R) \to \mathbb{R}$ is the energy functional given by
$$
J_R(u)=\frac 12 \int_{B_R}|\nabla u|^{2} +\frac 14\int_{B_R}\phi_u u^2dx-
\int_{B_R}F(u)dx,
$$
and
\[
\mathcal{M}_R=\{u\in H_0^1(B_R): J_R'(u)u^+=0=J_R'(u)u^-=0, u^\pm\neq
0\}.
\]
}
\begin{lemma} \label{ZZ1}
Let $c_0$ be the nodal level of $J$. Then
\[
\lim_{R \to +\infty}c_R=c_0.
\]
\end{lemma}

\noindent {\bf{Proof:}}
Since $R \mapsto c_R$ is a non-increasing function and {$c_R \geq c_0$ for all $R>0$}, if  $\lim_{R \to +\infty} c_R=\hat c>c_0$, then there exists $\varphi \in \mathcal M$ such that $J(\varphi)< \hat c$. From $\varphi \in \mathcal M$, $\varphi^\pm \neq 0$. Let $\varphi_n \in
C_0^\infty(\mathbb R^3)$ be such that $\varphi_n \to \varphi$ in $H^1(\mathbb
R^3)$. We may assume that $\varphi_n^\pm \neq 0$. By Lemma \ref{lema3}, there exist $t_n, s_n>0$ such that $J'(t_n\varphi_n^+ + s_n\varphi_n^-)\varphi_n^+ = 0$ and
$J'(t_n\varphi_n^+ + s_n\varphi_n^-)\varphi_n^- = 0$. In particular,
$
J'(t_n\varphi_n^+ + s_n\varphi_n^-)(t_n\varphi_n^+ + s_n\varphi_n^-) = 0.
$
Using that $(t_n\varphi_n^+ + s_n\varphi_n^-)^+ = t_n\varphi_n^+ \neq 0$ and $(t_n\varphi_n^+ + s_n\varphi_n^-)^- = s_n\varphi_n^- \neq 0$, we find that $t_n\varphi_n^+ + s_n\varphi_n^- \in \mathcal M\cap C_0^\infty(\mathbb R^3)$.  We claim that there exists a subsequence, still denoted by $(t_n\varphi_n^+ + s_n\varphi_n^-)$, such that $J(t_n\varphi_n^+ + s_n\varphi_n^-)\to J(\varphi)$.  Suppose for the moment that the limit holds. Let $n$ and $R>0$ such that $t_n\varphi_n^+ + s_n\varphi_n^- \in \mathcal  M_R$ and $J(t_n\varphi_n^+ + s_n\varphi_n^-) <
\hat{c}$. Hence,
\[
c_R \leq J(t_n\varphi_n^+ + s_n\varphi_n^-) < \hat c,
\]
and finally that
\[
\hat c = \lim_{R \to +\infty} c_R \leq J(t_n\varphi_n^+ + s_n\varphi_n^-) < \hat c,
\]
which is impossible. To establish the last claim, we start with the observation that there exist subsequences (not renamed) such that $t_n \to 1$ and $s_n \to 1$. In fact, suppose by contradiction that $\limsup_{n\to \infty} t_n > 1$. Given $\delta >0$ there exists a subsequence, still denoted by $t_n$, such that $t_n \geq \sigma$ for every $n$, for some $\sigma >1$. Since $J'(\varphi_n) \to J'(\varphi) = 0$ and the function $u\mapsto u^+$ is continuous, we have
\begin{equation}\label{eq3.1}
\|\varphi_n^+\|^2 + \int_{\mathbb{R}^3}\phi_{\varphi_n^+}(\varphi_n^+)^2dx  \leq \int_{\mathbb{R}^3}f(\varphi_n^+)\varphi_n^+ dx + o_n(1).
\end{equation}
On the other hand, $J'(t_n\varphi_n^+ + s_n\varphi_n^-)t_n\varphi_n^+ = 0$, that is
\begin{equation}\label{eq3.2}
\frac{1}{t_n^2}\|\varphi_n^+\|^2 + \int_{\mathbb{R}^3}\phi_{\varphi_n^+}(\varphi_n^+)^2dx = \int_{\mathbb{R}^3}\frac{f(t_n\varphi_n^+)t_n\varphi_n^+}{t_n^4} dx.
\end{equation}
Combining (\ref{eq3.1}) with (\ref{eq3.2}), gives
\begin{equation}\label{eq3.3}
 \left(1 - \frac{1}{t_n^2}\right)\|\varphi_n^+\|^2  {\leq} \int_{\mathbb{R}^3}\left[ \frac{f(\varphi_n^+)\varphi_n^+}{(\varphi_n^+)^4} - \frac{f(t_n\varphi_n^+)t_n\varphi_n^+}{(t_n\varphi_n^+)^4}\right](\varphi_n^+)^4 dx + o_n(1).
\end{equation}
From $(f_4)$ and Fatou's lemma, we have
\[
0 \leq \int_{\mathbb{R}^3}\left[ \frac{f(\sigma\varphi^+)\sigma\varphi^+}{(\sigma\varphi^+)^4} - \frac{f(\varphi^+)\varphi^+}{(\varphi^+)^4}\right](\varphi^+)^4 dx \leq
\left(\frac{1}{\sigma^2} - 1\right)\|\varphi^+\|^2 <0,
\]
which is impossible. Hence, $\displaystyle \limsup_{n\to \infty} t_n \leq 1$. Consequently, there exists a subsequence (not renamed) such that $\lim_{n\to \infty}t_n = t_0$. Taking to the limit as  $n\to \infty$ in  (\ref{eq3.3}) and using $(f_4)$ again, we get $t_0 =1$. In an exactly similar way,  there exists a subsequence (not renamed) such that $\displaystyle \lim_{n\to \infty}s_n = 1$. Finally, considering that
\begin{eqnarray*}
J(t_n\varphi_n^+ + s_n\varphi_n^-) &=& \frac{t_n^2}{2}\|\varphi_n^+\|^2 + \frac{s_n^2}{2}\|\varphi_n^-\|^2
 + \frac{t_n^4}{4}\int_{\mathbb{R}^3}\phi_{\varphi_n^+}(\varphi_n^+)^2dx\\ & & +
\frac{s_n^4}{4}\int_{\mathbb{R}^3}\phi_{\varphi_n^-}(\varphi_n^-)^2dx -  \int_{\mathbb{R}^3}F(t_n\varphi_n^+ + s_n\varphi_n^-)dx,
\end{eqnarray*}
we obtain that $J(t_n\varphi_n^+ + s_n\varphi_n^-)\to J(\varphi)$, by Lemma \ref{lm1} and the convergence $\varphi_n \to \varphi$ in $H^1(\mathbb
R^3)$. \qed


\section { The minimum level is achieved on $\mathcal M$  }

In this section, our main goal is to prove that the infimum $c_0$ of
$J$ on $\mathcal M$ is achieved. From Lemma \ref{lema2}(i), we
deduce that $c_0>0$. We start with this following lemma.

\begin{lemma}\label{lema4}
Suppose that $(u_n)$ be a sequence in $\mathcal M$ such that
$$
\limsup_{n\to \infty} J(u_n)<c+c_\infty.
$$
Then $u_n$ has a subsequence which converges weakly to some $w\in
H^1(\mathbb R^3)$ such that $w^\pm\neq 0$.
\end{lemma}

\noindent {\bf Proof.} From Lemma \ref{lema2}(i), $(u_n)$ is a
bounded sequence. Hence, without loss of generality, we can suppose
that there is $w \in H^1(\mathbb R^3)$ verifying
$u_n \rightharpoonup w$ in $H^{1}(\mathbb R^3)$ and
$u_n(x) \to w(x)$ almost everywhere in $\mathbb R^3$.  Observing that
\[
J(u_n)=J(u_n^+)+J(u_n^-)+\frac 12 \int_{\mathbb R^3}
\phi_{u^-_n}(u^+_n)^2dx,
\]
and
\[
J'(u_n^+)u^+_n=-\int_{\mathbb R^3}\phi_{u^-_n}(u^+_n)^2dx=
J'(u_n^-)u^-_n,
\]
we can suppose that
\[
J(u_n^+)+\frac 14 \int_{\mathbb
R^3}\phi_{u^-_n}(u^+_n)^2dx=\theta+o_n(1)
\]
and
\[
J(u_n^-)+\frac 14 \int_{\mathbb
R^3}\phi_{u^-_n}(u^+_n)^2dx=\sigma+o_n(1),
\]
where $\theta+\sigma<c+c_\infty$. We claim that $w^+ \neq 0$. Suppose by contradiction that $w^+\equiv 0$. From condition $(V_2)$ and Sobolev compact imbedding, we
have
$$
 \int_{\mathbb
R^3} V(x)(u^+_n)^2dx=\int_{\mathbb R^3} V_\infty(u^+_n)^2dx+o_n(1),
$$
which implies $J_\infty(u^+_n)= J(u^+_n)+o_n(1)$ and
$J_\infty'(u^+_n)u^+_n= J'(u^+_n)u^+_n+o_n(1)$. Hence,
\[
J_\infty(u_n^+)+\frac 14 \int_{\mathbb
R^3}\phi_{u^-_n}(u^+_n)^2dx=\theta+o_n(1)
\]
and
\[
J'_\infty(u_n^+)u^+_n=-\int_{\mathbb
R^3}\phi_{u^-_n}(u^+_n)^2dx+o_n(1).
\]
We observe that $\theta\geq  c_\infty$. In fact, let
$t_n>0$ be such that $J_\infty(t_nu^+_n)\geq J_\infty(tu^+_n)$, for
all $t>0$. We have three possibilities for  $(t_n)$:
\begin{enumerate}
\item [(i)] $\displaystyle \limsup_{n\to \infty}t_n>1$,
\item [(ii)]$\displaystyle \limsup_{n\to \infty} t_n=1$,
\item [(iii)]$\displaystyle \limsup_{n\to \infty}t_n<1$.
\end{enumerate}
We show now that (i) can not happen and (ii) or (iii) imply
$\theta\geq c_\infty$. From $J'_\infty(t_nu_n^+)t_nu^+_n=0$ we have
\begin{equation}\label{xy01}
t_n^2||u^+_n||_\infty^2+t_n^4\int_{\mathbb R^3}
\phi_{u^+_n}(u^+_n)^2dx=\int_{\mathbb R^3} f(t_nu^+_n)t_nu^+_ndx
\end{equation}
and from $J'(u_n)u^+_n=0$ follows
$$
||u^+_n||^2+\int_{\mathbb R^3} \phi_{u^+_n}(u^+_n)^2dx+\int_{\mathbb
R^3} \phi_{u^-_n}(u^+_n)^2dx=\int_{\mathbb R^3} f(u^+_n)u^+_ndx
$$
which implies
\begin{equation}\label{xy02}
||u^+_n||_\infty^2+\int_{\mathbb R^3}
\phi_{u^+_n}(u^+_n)^2dx+\int_{\mathbb R^3}
\phi_{u^-_n}(u^+_n)^2dx=\int_{\mathbb R^3} f(u^+_n)u^+_ndx+o_n(1).
\end{equation}
Combining (\ref{xy01}) and (\ref{xy02}) we get
\begin{equation}\label{xy03}
(1-\frac 1{t_n^2})||u^+_n||_\infty^2+\int_{\mathbb R^3}
\phi_{u^-_n}(u^+_n)^2dx=\int_{\mathbb R^3}\left [
\frac{f(u^+_n)}{(u^+_n)^3}-\frac{f(t_nu^+_n)}{(t_nu^+_n)^3}\right
](u^+_n)^4dx+o_n(1).
\end{equation}
If (i) holds, there exists  $a>1$ such that $t_n\geq a$ for
infinitely many $n$. By Lemma \ref{lema2} (ii), the left hand
in (\ref{xy03}) is bounded from below by a positive number. On the other hand, by $(f_4)$,
the integral in the the right hand of (\ref{xy03}) is
non-positive. This yields a contradiction. Hence (i) does not hold.
Suppose that (iii) holds. Then, $t_n\leq1$ and Remark \ref{rem} imply
\begin{eqnarray*}
4c_\infty &\leq& 4J_\infty(t_nu^+_n)=4J_\infty(t_nu^+_n)-J'_\infty(t_nu^+_n)(t_nu^+_n)\\
&=& t_n^2{||u^+_n||_*^2}+\int_{\mathbb R^3}
[f(t_nu^+_n)t_nu^+_n-4F(t_nu^+_n)]dx\\
&\leq& {||u^+_n||_*^2}+\int_{\mathbb R^3} [f(u^+_n)u^+_n-4F(u^+_n)]dx\\
&=& 4J_\infty(u^+_n)-J'_\infty(u^+_n)(u^+_n)\\
&=& 4J_\infty(u^+_n)+\int_{\mathbb R^3}\phi_{u^-_n}(u^+_n)^2dx+o_n(1)\\
&=&4\theta+o_n(1).
\end{eqnarray*}
Taking to the limit $n\to \infty$, we find $\theta\geq c_\infty$.
If (ii) occurs, there exists a subsequence (still denoted by $t_n$) such that $\lim_{n\to \infty}t_n =1$. As a consequence
\[
4J_\infty(t_nu_n^+) - J_\infty'(t_nu_n^+)(t_nu_n^+) = 4J_\infty(u_n^+) - J_\infty'(u_n^+)(u_n^+) + o_n(1).
\]
Thus,
\begin{eqnarray*}
4c_\infty &\leq&
4J_\infty(t_nu^+_n)=4J_\infty(t_nu^+_n)-J'_\infty(t_nu^+_n)(t_nu^+_n)\\
&=&  4J_\infty(u_n^+) - J_\infty'(u_n^+)(u_n^+) + o_n(1)\\
&=& 4J_\infty(u^+_n)+\int_{\mathbb
R^3}\phi_{u^-_n}(u^+_n)^2dx+o_n(1)\\
&=& 4\theta+o_n(1).
\end{eqnarray*}
Taking the limit $n\to \infty$, we also obtain $\theta \geq c_\infty$.
Since $\theta+\sigma<c+c_\infty$ and  $\theta \geq c_\infty$, we have  $\sigma<c$. Let $s_n>0$ be such that $J(s_nu^-_n)\geq J(tu^-_n)$,
for all $t>0$. Using that $J'(u^-_n)(u^-_n)<0$, we get $s_n<1$. Hence,
\begin{eqnarray*}
4c&\leq&
4J(s_nu^-_n)=4J(s_nu^-_n)-J'(s_nu^-_n)(s_nu^-_n)\\
 &=& s_n^2\|u^-_n\|^2+\int_{\mathbb R^3} [f(s_nu^-_n)s_nu^-_n-4F(s_nu^-_n)]dx\\
 &\leq& \|u^-_n\|^2+\int_{\mathbb R^3} [f(u^-_n)u^-_n-4F(u^-_n)]dx\\
 &=&4J(u^-_n)-J'(u^-_n)(u^-_n)\\
 &=&4J(u^-_n)+\int_{\mathbb R^3}\phi_{u^-_n}(u^+_n)^2dx\\
 &=&4\sigma+o_n(1),
\end{eqnarray*}
which implies  $c \leq \sigma$, contrary to  $\sigma < c$. Hence, $w^+\neq 0$ as claimed. Similar arguments to those above show that
 $w^-\neq 0$, and the proof is complete.\qed

\begin{lemma}\label{lema5}
If $c_0<c+c_\infty$, there exists a $w\in \mathcal M$ which
minimizes $J$ on $\mathcal M$.
\end{lemma}

 \noindent {\bf Proof.}  By Proposition \ref{appendix}, there exists a least energy sign-changing solution $u_n$ to  $(AP_R)$, for $R=n$, that is
 $ J(u_n) = c_n= \inf_{ \mathcal M_n}J$, where  $c_n=c_R$ and $\mathcal M_n=\mathcal M_R$. By Lemma \ref{ZZ1},   $c_n \to c_0$ as $n\to \infty$. Moreover,
  $J'(u_n)v=0$ for all $v \in H_0^1(B_n)$. Since  $c_0 < c + c_\infty$, $u_n$  converges weakly to some $w\in H^1(\mathbb R^3)$ such that
$w^\pm\neq 0$, by Lemma \ref{lema4}. Using that $J'(u_n)v=0$ for all $v \in H_0^1(B_n)$, we get $J'(w)=0$ and consequently $w\in \mathcal M$.
 We claim that $J(w)=c_0$. In fact, combining Fatou's lemma with Remark \ref{rem}, we have
$$
c_0 \leq  J(w)-\frac{1}{4}J'(w).w\leq \liminf_{n\to \infty } \left (J(u_n)
-\frac{1}{4}J'(u_n)u_n \right)=\lim_{n\to \infty} J(u_n)=c_0,
$$
which implies that  $c_0 = J(w)$. \qed

Until this moment, we have proved that under condition
$c_0<c+c_\infty$, there exists a $w\in \mathcal M$, such that
$J(w)=c_0$ and $J'(w)=0$.

\section{Estimate on the level $c_0$}

This section is devoted to show that $c_0<c+c_\infty$.  The proofs
in this section are based upon ideas found in \cite{lww}. From now
on fix $w, v \in H^1(\mathbb R^3)$ \emph{ground state solution} of
$(P)$ and $(P_\infty)$  given by \cite[Theorems 1.5 and 1.3]{JMAA2010}
respectively.   We know that $w$ and $v$ should have defined sign.
Without loss of generality we will suppose that:
$$
w>0,\,\, v>0, \mbox { in }\mathbb R^3, \, J(w)=c,
\,J_\infty(v)=c_\infty, \, J'(w)=0, \mbox { and }J'_\infty(v)=0.
$$

Using Moser's iteration and De Giorgi's iteration, we can show that  $w$ and $v$  have exponential decay, and consequently,
$\phi_v$ and $\phi_w$ have the same behavior. More precisely:
\begin{lemma}\label{lm8}
There exist $C>0$ and $\delta>0$ such that for all $R>0$:
\begin{eqnarray*}
&&\int_{|x|\geq R} (|\nabla w|^2+w^2)dx \leq Ce^{-\delta R},\,
\int_{|x|\geq R} (|\nabla v|^2+v^2)dx \leq Ce^{-\delta R}, \\
&&\int_{|x|\geq R} (F(w)+wf(w)+F(v)+v f(v))dx \leq
Ce^{-\delta R},\\
&&\int_{|x|\geq R} \phi_wv^2dx+\int_{|x|\geq R} \phi_v w^2dx \leq
Ce^{-\delta R}.
\end{eqnarray*}
\end{lemma}

\medskip
For each $n\in \mathbb N$, set $v_n(x)=v(x+ne_1)$, where
$e_1=(1,0,0)\in \mathbb R^3$. The same conclusion of Lemma
\ref{lm8} is satisfied by  function $v_n$ and
\begin{equation}\label{nltexp}
\int_{\mathbb R^3} \phi_wv_n^2dx=\int_{\mathbb R^3}
\phi_{v_n}w^2dx=O(e^{-n\delta }).
\end{equation}

\begin{lemma}\label{lm9}
Suppose that $V$ satisfies $(V_2)-(V_3)$ and $f$ satisfies $(f_2)$
and $(f_5)$. Then,
\[
\sup_{(\alpha,\beta)\in \mathbb R^2}J(\alpha w + \beta
v_n)<c+c_\infty,
\]
provided  $n$ is sufficiently large.
\end{lemma}

\noindent {\bf{Proof:}} We start proving that there is $r_o>0$ such
that $J(\alpha w + \beta v_n)\leq 0$, for all $(\alpha,\beta)\in
\mathbb R^2$ such that $\alpha^2+\beta^2\geq r_o$, and  $n\geq r_o$.
Since $J(v)\leq J_\infty(v) $, for all $v$, it is sufficient to show that $J_\infty(\alpha w + \beta
v_n)\leq 0$, for all $\alpha^2+\beta^2\geq r_o$, $n\geq r_o$.  In fact, suppose that $J_\infty$ does not
satisfies this claim. Thus, for each $n$, there are $(\alpha_n,\beta_n)\in
\mathbb R^2$ such that $J_\infty(\alpha_n w + \beta_n v_n)> 0$ and
$\alpha_n^2+\beta_n^2\to \infty$, that is,
\begin{equation}\label{eq1113}
\frac 12 {||\alpha_n w + \beta_n v_n||_*^2} +\frac
14\int_{\mathbb R^3}\phi_{(\alpha_n w + \beta_n v_n)} (\alpha_n w +
\beta_n v_n)^2dx\geq \int_{\mathbb R^3} F(\alpha_n w + \beta_n
v_n)dx,
\end{equation}
We have {$||v_n||_*=||v||_*$}, and from Lemma \ref{lm8}
\begin{equation}\label{o1}
\int_{\mathbb R^3} (\nabla w \nabla v_n+V_\infty
wv_n)dx=O(e^{-n\delta}).
\end{equation}
It follows that
\begin{equation}\label{o2}
{||\alpha_n w + \beta_n
v_n||_*^2=\alpha_n^2||w||_*^2+\beta_n^2||v||_*^2+O(e^{-n\delta})}
\end{equation}
and then {$\sigma_n=||\alpha_n w + \beta_n v_n||_*\to +\infty$.}
Set
\[
{z_n=\frac{\alpha_n w + \beta_n v_n}{||\alpha_n w + \beta_n
v_n||_*},}
\]
and suppose that $z_n \rightharpoonup z$. Dividing (\ref{eq1113}) by
$\sigma_n^4$,  we have
\[
\frac 1{2\sigma_n^2}+\frac14\int_{\mathbb R^3} \phi_{z_n} z_n^2dx
\geq \int_{\mathbb R^3} \frac{F(\alpha_n w + \beta_n v_n)}{(\alpha_n
w + \beta_n v_n)^4}z_n^4dx.
\]
{the boundedness of $(z_n)$ together with the above inequality and $(f_3)$ }shows that $z\equiv 0$. Passing to the
limit as $n\to \infty$ in the equality below
{\begin{eqnarray*}
o_n(1)=\lefteqn{\int_{\mathbb R^3} (\nabla w \nabla z_n+V_\infty wz_n)dx=
\alpha_n||\alpha_n w  +
  \beta_n v_n||_*^{-1}||w||_*^2  }\\    && +\   \beta_n||\alpha_n w + \beta_n v_n||_*^{-1}\int_{\mathbb R^3}
   (\nabla w \nabla v_n+V_\infty wv_n)dx,
\end{eqnarray*}}
we obtain from (\ref{o1}) and (\ref{o2}) that {$\alpha_n||\alpha_n w
+ \beta_n v_n||_*^{-1}$} converges to $0$. By Lemma
\ref{lm1}$ (iv)$,  $J_\infty(\alpha_n w + \beta_n
v_n)=J_\infty(\alpha_n w_n + \beta_n v)$ where $w_n(x)=w(x-ne_1)$.
Proceeding exactly as the previous argument, we can show that
{$\beta_n||\alpha_n w + \beta_n v_n||_*^{-1}$} converges to $0$.
From (\ref{o2}), $z_n \to 0$, which contradicts {$||z_n||_*=1$}.
Hence, the claim holds for $J_\infty$, and, in consequence, for $J$.

Now we consider $n\geq r_o$, $\alpha^2+\beta^2\leq r_o$. From Lemma
\ref{lm8},  there are $\delta>0$ and $C=(w,v,r_o)$ such that
\begin{equation}\label{eq1114}
J(\alpha w+\beta v_n)\leq J(\alpha w)+J(\beta
v_n)+Ce^{-n\delta}.
\end{equation}
In fact, from (\ref{nltexp}) and Lemma
\ref{lm8}, we have
\begin{eqnarray*}
|\int_{\mathbb R^3} F(\alpha w+\beta v_n)dx-\int_{\mathbb R^3}
F(\alpha w)dx-\int_{\mathbb R^3} F(\beta v_n)dx|\leq
Ce^{-n\delta}, \\
|\int_{\mathbb R^3} \phi_{(\alpha w+\beta v_n)}(\alpha
w+\beta v_n)^2dx-\alpha^4\int_{\mathbb R^3}
\phi_ww^2dx-\beta^4\int_{\mathbb R^3}
\phi_{v_n}v_n^2dx|\leq Ce^{-n\delta}
\end{eqnarray*}
and $\left |||\alpha w+\beta v_n||^2- \alpha^2 ||w||^2-\beta^2
||v_n||^2\right |\leq Ce^{-n\delta}.$ The condition $(V_3)$
together with Lemma \ref{lm8} imply that
\begin{equation} \label{eq1115}
J(\beta v_n)\leq J(t_nv_n)\leq J_\infty(t_n
v_n)+Ce^{-n\delta}-\rho(n+1)t_n^2|v|_{L^2(B_1(0))}^2,
\end{equation}
where $t_n>0$ is such that $J(t_nv_n)\geq J(tv_n)$, for
all $t>0$ and $\rho $ is the non-decreasing function given by
$(V_3)$. In order to verify (\ref{eq1115}), we first observe that
\[
J(t_n v_n)= J_\infty(t_n v_n)+\int_{\mathbb R^3}
(V(x)-V_\infty)t_n^2v_n^2dx
\]
and
\[
\int_{\mathbb R^3}
(V(x)-V_\infty)t_n^2v_n^2dx\leq\int_{|x-ne_1|\leq 1}
(V(x)-V_\infty)t_n^2v_n^2dx.
\]
For $R_o$  given by $(V_3)$, set $n \geq R_o+1$. If $|x-ne_1|\leq
1$, we have $n-1\leq |x|\leq n+1$ and then $ -\rho(|x|)\leq
-\rho(n+1)$. Hence,
\[
\int_{|x-ne_1|\leq 1}(V(x)-V_\infty)t_n^2v_n^2dx \leq
-\rho(n+1)\int_{|x-ne_1|\leq 1}t_n^2v_n^2dx,
\]
and so (\ref{eq1115}) is justified.    By the definition of $t_n$,
we have
\begin{equation}\label{auxiliar}
t_n\|v_n\|^2 +t_n^3\int_{\mathbb R^3} \phi_{v_n}
v_n^2dx= \int_{\mathbb{R}^3} f(t_nv_n)v_ndx =
\int_{\mathbb{R}^3} f(t_nv)v dx.
\end{equation}
Combining (\ref{auxiliar}) with the fact {$\|v_n\|^2=\|v_n\|^2_* +
o_n(1)= \|v\|^2_* + o_n(1)$}, by $(V_2)$-$(V_3)$, and using Lemma
\ref{lm1}(iv) and (\ref{estimativa}), we get
\begin{eqnarray*}
 t_n^2(\|v\|^2_* +
o_n(1)) &\leq&  t_n^2(\|v\|^2_* + o_n(1))+t_n^4\int_{\mathbb
R^3} \phi_{v} v^2dx = t_n\int_{\mathbb{R}^3} f(t_nv)v dx\\
&\leq & \epsilon t_n^2\int_{\mathbb {R}^3}v^{2} dx +  Ct_n^6\int_{\mathbb {R}^3}v^{6} dx,
\end{eqnarray*}
for some positive constant $C$. Therefore, there exists $\tau > 0$
such that  $t_n^2 \geq \tau$ for every $n$. Using (\ref{eq1114}),  (\ref{eq1115}) and the fact that $J_\infty(t_nv_n) = J_\infty(t_nv) \leq J_\infty(v) = c_\infty$,  we have
\begin{eqnarray*}
J(\alpha w+\beta v_n)&\leq& J(\alpha w)+J_\infty(t_n v_n)+Ce^{-n\delta}-t_n^2\rho(n+1)|v|_{L^2(\mathbb R^3)}^2\\
&\leq& c+c_\infty+e^{-n\delta}\left(C-\tau |v|_{L^2(\mathbb
R^3)}^2 e^{n\delta}\rho(n+1)\right)
\end{eqnarray*}
and the proof follows by the limit condition on $\rho$ in $(V_3)$.
\qed

 \medskip

We have the following lemma:
\begin{lemma}{\label{c*<c}} The number $c_0$ verifies the following inequality
\begin{equation}
c_0 < c+c_\infty.
\end{equation}
\end{lemma}

\noindent {\bf{Proof:}} Let  $w$ and $v_n$ be functions as in
the proof of Lemma \ref{lm9}. Let $D= [1/2,3/2]\times [1/2,3/2]$ and
\[
\Psi(\xi,\tau) = \left(J'((\xi w - \tau v_n)^+)(\xi w - \tau v_n)^+, J'((\xi w - \tau v_n)^-)(\xi w - \tau v_n)^-\right).
\]
Using that $J'(w)w=0$ and $(f_4)$, we obtain
\begin{equation}\label{m1}
J'(\frac{1}{2}w)\frac{1}{2}w>0\quad \mbox{and}\quad  J'(\frac{3}{2}w)\frac{3}{2}w < 0.
\end{equation}
Property (iv) of Lemma \ref{lm1}, condition $(V_2)$ and the fact that $J'_\infty(v)v=0$ imply that there exists $n_0 \in \mathbb{N}$ such that
\begin{equation}\label{m2}
J'(\frac{1}{2}v_n)\frac{1}{2}v_n>0\quad \mbox{and}\quad  J'(\frac{3}{2}v_n)\frac{3}{2}v_n < 0,\quad \forall\, n\geq n_0.
\end{equation}
Since  $v(x) \to 0$ as $|x| \to \infty$, it follows from (\ref{m1})-(\ref{m2}),  by increasing $n_0$ if necessary,  that
\begin{equation}\label{m3}
J'((\frac{1}{2}w - \tau v_n)^+)(\frac{1}{2}w - \tau v_n)^+ >0 \mbox{ and }
J'((\frac{3}{2}w - \tau v_n)^+)(\frac{3}{2}w - \tau v_n)^+ <0,
\end{equation}
for every $n \geq n_0$ and $\tau \in [1/2,3/2]$, and
\begin{equation}\label{m4}
J'((\xi w - \frac{1}{2} v_n)^-)(\xi w - \frac{1}{2}v_n)^- >0 \mbox{ and }
J'((\xi w - \frac{3}{2} v_n)^-)(\xi w - \frac{3}{2}v_n)^- <0,
\end{equation}
for every $n \geq n_0$ and $\xi \in [1/2,3/2]$. Noting that the function $\Psi$ is continuous in $D$ and considering the inequalities (\ref{m3})-(\ref{m4}), we can apply Miranda's theorem \cite{Miranda} and conclude that there exists $(\xi_0, \tau_0) \in D$ such that $\Psi(\xi_0, \tau_0) = (0,0)$. This gives $\xi_0 w - \tau_0 v_n \in \mathcal{M}$ for every $n \geq n_0$. Consequently,
$$
c_0 \leq J(\xi_0 w - \tau_0 v_n ),
$$
which  implies
$$
c_0 \leq \sup_{(\alpha,\beta)\in \mathbb R^2}J(\alpha w + \beta
\omega_n).
$$
The lemma follows combining the last inequality with  Lemma
\ref{lm9}. \qed

 \section{Theorem \ref{fth}}

In this section we establish a proof of Theorem \ref{fth}. From Sections 5 and 6, there exists a critical point  $w$ of $J$, which is a sign-changing solution for
problem $(SP)$. The proof is completed by showing that  $w$ has exactly two nodal
domains. Arguing by contradiction, we suppose  that
$
w=u_1+u_2+u_3,
$
with
$
u_i \not= 0, u_1 \geq 0, u_2 \leq 0$  and ${\rm{supp}}(u_i) \cap {\rm{supp}}(u_j) = \emptyset$ for $i \not= j, i,j=1,2,3$, with ${\rm{supp}}(u_i)$ denoting the support of $u_i$. Setting $v = u_1 + u_2$, we see that $v^{\pm} \not= 0$. Moreover,
using the fact that $J'(w)=0 $, it follows that
$$
J'(v)(v^{\pm}) \leq 0.
$$
By Lemma \ref{lema3}, there are $t,s \in (0,1]$ such that
$
tv^+ + sv^{-} \in \mathcal{M}
$
or equivalently,
$
tu_1 + su_2 \in \mathcal{M},
$
and so,
\begin{equation}\label{eq30}
J(tu_1+su_2) \geq c_0.
\end{equation}
Since $w = v + u_3$, we have $w^2 = v^2 + u_3^2$ and $\phi_w = \phi_v + \phi_{u_3}$. Hence,
\begin{equation}\label{eq31}
J(w) = J(v) + J(u_3) + \frac{1}{2}\int_{\mathbb{R}^3}\phi_vu_3^2 dx.
\end{equation}
Supposing that $u_3\neq 0$, we claim that
\begin{equation}\label{eq32}
J(u_3) + \frac{1}{4}\int_{\mathbb{R}^3}\phi_vu_3^2 dx > 0.
\end{equation}
 In fact, by Remark \ref{rem} and using $J'(w)u_3 = 0$ combined with $u_3 \not=0$, we obtain
\begin{eqnarray*}
J(u_3) + \frac{1}{4}\int_{\mathbb{R}^3}\phi_vu_3^2 dx &=& J(u_3) + \frac{1}{4}\int_{\mathbb{R}^3}\phi_vu_3^2 dx - \frac{1}{4}J'(w)u_3\\
 &=& \frac{1}{4}\|u_3\|^2 + \frac{1}{4}\int_{\mathbb{R}^3} (f(u_3)u_3 - 4F(u_3))dx >0.
\end{eqnarray*}
Similar arguments to those above show that
\begin{equation}\label{eq33}
J(v) + \frac{1}{4}\int_{\mathbb{R}^3}\phi_vu_3^2 dx = \frac{1}{4}\|v\|^2 + \frac{1}{4}\int_{\mathbb{R}^3} (f(v)v - 4F(v))dx.
\end{equation}
From (\ref{eq30})-(\ref{eq33}), for every $t, s \in (0,1]$, we have
\begin{eqnarray*}
c_0 &\leq & J(tu_1+su_2) = J(tu_1+su_2)  - \frac{1}{4}J'(tu_1+su_2)(tu_1+su_2)\\
&=& \frac{t^2}{4}\|u_1\|^2 + \frac{s^2}{4}\|u_2\|^2 + \frac{1}{4}\int_{\mathbb{R}^3}(f(tu_1 + su_2)(tu_1 + su_2) - 4F(tu_1 + su_2))dx\\
&\leq &\frac{1}{4}\|u_1\|^2 + \frac{1}{4}\|u_2\|^2 + \frac{1}{4}\int_{\mathbb{R}^3}(f(u_1 + u_2)(u_1 + u_2) - 4F(u_1 + u_2))dx\\
&=& J(v) + \frac{1}{4}\int_{\mathbb{R}^3}\phi_vu_3^2 dx \\
&<& J(v) + \frac{1}{4}\int_{\mathbb{R}^3}\phi_vu_3^2 dx + J(u_3) + \frac{1}{4}\int_{\mathbb{R}^3}\phi_vu_3^2 dx\\
&=& J(v) + J(u_3)  + \frac{1}{2}\int_{\mathbb{R}^3}\phi_vu_3^2 dx\\
&=& J(w)\\
&=& c_0,
\end{eqnarray*}
which is a contradiction. Therefore,  $u_3=0$ and $w$ has exactly
two nodal domains. \qed

\section{Appendix}

In this appendix we present the existence of a least energy sign-changing solution for the following system:
$$
\left\{\begin{array}{l}
- \Delta u +V(x)u+\phi u= f(u),  \mbox{ in }  B_R, \\
-\Delta \phi=\tilde{u}^2,\mbox{ in }  \mathbb R^3,\\
 \phi\in D^{1,2}(\mathbb R^3),\  u \in H_0^1(B_R),
\end{array}\right. \leqno (AP_R)
$$
where $B_R$ is the of radius $R$ centered at $0$ and
$$
\tilde{u}(x)=
\left\{
\begin{array}{cl}
u(x) & \mbox{if $x\in B_R$}, \\
0 & \mbox{if $x\in R^3\setminus B_R$}.
\end{array}
\right.
$$

The energy functional $J_R:H^{1}_{0}(B_R) \to \mathbb{R}$ associated with $(AP_R)$ is given by
$$
J_R(u)=\frac 12 \int_{B_R}|\nabla u|^{2} +\frac 14\int_{B_R}\phi_u u^2dx-
\int_{B_R}F(u)dx,
$$
In this section, we will use $J$ to denote $J_R$

It is important to point out that the problem $(AP_R)$ is not the
same problem considered in \cite{nodal1}, namely
$$
\left\{\begin{array}{l}
- \Delta u +V(x)u+\phi u= f(u),  \mbox{ in }  B_R, \\
-\Delta \phi= u^2,\mbox{ in }   B_R,\\
 \phi, u \in H_0^1(B_R).
\end{array}\right. \leqno (P')
$$
 The local terms are different, and
consequently, the associated functionals are different. Here
$\phi_u$ is the restriction of $\phi_{\tilde{u}}$  to $B_R$ and, for $(P')$,
$\phi_u$ should be a function in the Sobolev space $H^1_0(B_R)$.

\begin{proposition}\label{appendix}
Suppose that $f$ satisfies $(f_1)-(f_4)$. Then, for any $R>0$, problem $(AP_R)$ possesses a least energy sign-changing solution.
\end{proposition}

\noindent {\bf{Proof:}} This proof follows \cite{nodal1}. Let $(w_n)$ be a sequence in $\mathcal M_R$ such that
$$
\lim_{n\to \infty} J(w_n)=c_R.
$$
Lemma \ref{lema2}(i) shows that $(w_n)$ is a bounded sequence. Hence, without loss of generality, we can suppose that there is $w
\in H_0^1(B_R)$ verifying
$$
w_n \rightharpoonup w \,\,\, \mbox{in} \,\,\, H_0^{1}(B_R),
$$
and
$$
w_n(x) \to w(x) \,\,\, \mbox{a.e. in} \,\, B_R.
$$
The condition $(f_2)$ combined with the {compactness lemma of
Strauss} \cite[Theorem A.I, p.338]{bl} gives
$$
\lim_{n\to \infty} \int_{B_R}|w_n^{\pm}|^{p}dx = \int_{B_R}|w^{\pm}|^{p}dx,
$$
$$
\lim_{n\to \infty} \int_{B_R} w_n^\pm f(w_n^\pm)dx=\int_{B_R} w^\pm f(w^\pm)dx,
$$
and
$$
\lim_{n\to \infty} \int_{B_R}  F(w_n^\pm)dx=\int_{B_R} F(w^\pm)dx.
$$
From $(f_1)-(f_2)$, given $\epsilon >0$ and $q \in (2, 6)$, there exists $C=C(\epsilon, q)$ such that
\[
f(s)s \leq \epsilon(|s|^2 + |s|^6) + C|s|^q, \forall s\in \mathbb{R}.
\]
As $w_n \in \mathcal M_R$, we have $J'(w_n)w_n^\pm <0$. Combining this fact with Lemma  \ref{lema2}, we get
\[
\rho^2 \leq \|w_n^{\pm}\|^2 < \int_{B_R}f(w_n^{\pm})w_n^{\pm}dx \leq \epsilon\int_{B_R}(|w_n^{\pm}|^2 + |w_n^{\pm}|^6)dx + C\int_{B_R}|w_n^{\pm}|^q dx.
\]
Using that $(w_n)$ is bounded in $H^1_0(B_R)$, it follows from the Sobolev imbeddings that there exists $C_1 >0$ such that
\[
\rho^2 \leq \epsilon C_1 + C\int_{B_R}|w_n^{\pm}|^q dx.
\]
For fixed $\epsilon = \rho^2/2C_1$, we find
\[
\int_{B_R}|w_n^{\pm}|^q dx \geq \frac{\rho^2}{2C},
\]
which shows that
\[
\liminf_{n\to \infty} \int_{B_R}|w_n^{\pm}|^q dx >0,
\]
and consequently  $w^\pm\neq 0$. Then, by Lemma \ref{lema3} there are $t,s>0$
verifying
\begin{equation}\label{eq0}
J'(tw^++sw^-)w^+=0 \,\,\, \mbox{and} \,\,\, J'(tw^++sw^-)w^-=0.
\end{equation}
We claim that $t,s\leq 1$. In fact, since $J'(w_n)w_n^{\pm} =0$,
\begin{eqnarray*}
\|w_n^+\|^2 + \int_{{B_R}}\phi_{w_n^+}(w_n^+)^2dx + \int_{{B_R}}\phi_{w_n^-}(w_n^+)^2dx &=& \int_{{B_R}}f(w_n^+)w_n^+dx\\
\|w_n^-\|^2 + \int_{{B_R}}\phi_{w_n^-}(w_n^-)^2dx + \int_{{B_R}}\phi_{w_n^+}(w_n^-)^2dx &=& \int_{{B_R}}f(w_n^-)w_n^-dx.
\end{eqnarray*}
Taking $n\to \infty$, we obtain
\begin{equation}
\|w^+\|^2 + \int_{{B_R}}\phi_{w^+}(w^+)^2dx + \int_{{B_R}}\phi_{w^-}(w^+)^2dx \leq \int_{{B_R}}f(w^+)w^+dx \label{ap1}
\end{equation}
\begin{equation}
\|w^-\|^2 + \int_{{B_R}}\phi_{w^-}(w^-)^2dx + \int_{{B_R}}\phi_{w^+}(w^-)^2dx \leq \int_{{B_R}}f(w^-)w^-dx \label{ap2}.
\end{equation}
From (\ref{eq0}),
\begin{equation}
t^2\|w^+\|^2 + t^4\int_{{B_R}}\phi_{w^+}(w^+)^2dx + t^2s^2\int_{{B_R}}\phi_{w^-}(w^+)^2dx = \int_{{B_R}}f(tw^+)tw^+dx \label{ap3}
\end{equation}
\begin{equation}
s^2\|w^-\|^2 + s^4\int_{{B_R}}\phi_{w^-}(w^-)^2dx + t^2s^2\int_{{B_R}}\phi_{w^+}(w^-)^2dx = \int_{{B_R}}f(sw^-)sw^-dx \label{ap4}.
\end{equation}
There is no loss of generality in assuming that $s\geq t$. From (\ref{ap4}), we have
\begin{equation}
s^2\|w^-\|^2 + s^4\int_{{B_R}}\phi_{w^-}(w^-)^2dx + s^4\int_{{B_R}}\phi_{w^+}(w^-)^2dx \geq \int_{{B_R}}f(sw^-)sw^-dx \label{ap5}.
\end{equation}
Combining (\ref{ap2}) with (\ref{ap5}), we find
\begin{equation}
\left(\frac{1}{s^2}-1\right)\|w^-\|^2 \geq \int_{{B_R}}\left( \frac{f(sw^-)sw^-}{(sw^-)^4} - \frac{f(w^-)w^-}{(w^-)^4}\right)(w^-)^4dx. \label{ap6}
\end{equation}
Whenever $s>1$, the left side of (\ref{ap6}) is negative, whereas the right side of (\ref{ap6}) is positive, by $(f_4)$. Hence, $s\leq 1$ and, in consequence, $t\leq 1$, as claimed. From (\ref{eq0}), $tw^+ + sw^- \in \mathcal{M}_R$, which gives
\begin{eqnarray}
c_R &\leq &  J( tw^+ + sw^- ) = J(tw^+ + sw^-) - \frac{1}{4}J'(tw^+ + sw^-)(tw^+ + sw^-)\nonumber\\
&=&( J(tw^+) - \frac{1}{4}J'(tw^+)(tw^+) ) + (J(sw^-) - \frac{1}{4}J'(sw^-)(sw^-))\nonumber\\
&\leq& [J(w^+) - \frac{1}{4}J'(w^+)(w^+)] + [J(w^-) - \frac{1}{4}J'(w^-)(w^-)] \nonumber \\
&=& J(w) - \frac{1}{4}J'(w)(w), \label{ap7}
\label{ap7}
\end{eqnarray}
where we have used Remark \ref{rem} and the above claim. Using Fatou's lemma and Remark \ref{rem} and (\ref{ap7}), we get
\[
c_R \leq J( tw^+ + sw^- ) \leq \liminf_{n\to \infty}\left[ J(w_n) - \frac{1}{4}J'(w_n)(w_n)\right] = c_R,
\]
yields $J(tw^++sw^-)=c_R$, that is, $w_R \doteq tw^++sw^-$ minimizes $J$ on $\mathcal M_R$. To conclude the proof, it only remains to verify that  $w_R$ is a critical point of $J$. Suppose, by contradiction, that $J'(w_R)\neq 0$. Thus, there exist $\alpha >0$ and $v_0\in H_0^1(B_R)$, $\|v_0\| =1$, such that
\[
J'(w_R)v_0 = 2\alpha >0.
\]
Since $J\in C^1(H_0^1(B_R),\mathbb{R})$ and the functions $k^\pm: H_0^1(B_R) \to H_0^1(B_R)$ defined by  $k^\pm(u) = u^\pm$ are continuous,  there exists $r>0$ such that
\[
J'(v)v_0 > \alpha,\   v^\pm \neq 0,\ \mbox{ for every $v\in H_0^1(B_R),\  \|v - w_R\| \leq r$.}
\]
Fix $0 < \xi < 1 < \chi$. Let $D= (\xi,\chi)\times (\xi,\chi) \subset \mathbb{R}^2$, $h(a,b)=J(aw_R^+ + bw_R^-)$  and $\Phi(a,b) = \left(\frac{\partial h}{\partial a}(a,b),  \frac{\partial h}{\partial b}(a,b)\right)$, for $(a,b)\in \overline D$.  By \cite[Lemma 2.4]{nodal1},  we have
\begin{itemize}
\item[(a)] $h(a,b) < h(1,1) = J(w_R)$, for every $a, b \geq 0$ such that $(a,b) \neq (1,1)$.
\item[(b)] ${\rm{det}}(\Phi')(1,1) >0$.
\end{itemize}
We can take $0 < \xi < 1 < \chi$ such that
\begin{itemize}
\item[(i)] $(1,1)\in D$ and $\Phi(a,b) = (0,0)$ if, and only if, $a=b=1$.
\item[(ii)] $c_R \not\in h(\partial D)$.
\item[(iii)] $\{ aw_R^+ + bw_R^-\, :\, (a,b) \in \overline {D}\} \subset B_r(w_R)$.
\end{itemize}
Since $J$ is continuous, there exists $r'>0$ such that $B=\overline{B_{r'}(w_R)}\subset B_r(w_R)$ and $B\cap\{aw_R^+ + bw_R^-\, :\,  (a,b) \in \partial D\} = \emptyset$. Let $\rho(u) = {\rm{dist}}(u, B^c)$, $u\in H_0^1(B_R)$. Set the bounded Lipschitz vector field $V(u) = -\rho(u)v_0$, $u\in H_0^1(B_R)$. For each $u\in H_0^1(B_R)$, let $\eta(\tau) = \eta(\tau, u)$ denote the solution of
\[
\eta'(\tau) = V(\eta(\tau)), \tau >0,\quad \eta(0) = u.
\]
Note that
\begin{enumerate}
\item If $u \not\in B$, $\eta(\tau, u) = u$, for every $\tau \geq 0$;
\item If $u \in B$, the function $\tau \to J(\eta(\tau, u))$ is decreasing and $\eta(\tau, u) \in B$, for every $\tau \geq 0$;
\item There exists $\tau_0 >0$ such that $J(\eta(\tau,w)) \leq J(w) - ((r'\alpha)/2)\tau$, for every $\tau \in [0, \tau_0]$.
\end{enumerate}
Define $\gamma(a,b) = \eta(\tau_0, aw_R^+ + bw_R^-)$, for  $(a,b)\in \overline D$. We have,
\[
J(\gamma(a,b)) \leq h(a,b) < h(1,1) = J(w_R) = c_R,\ \forall\, (a,b) \in \overline{D}\setminus\{(1,1)\},
\]
and
\[
J(\gamma(1,1)) = J(\eta(\tau_0, w_R^+ + w_R^-)) = J(\eta(\tau_0, w_R)) < J(\eta(0, w_R))= J(w_R) = c_R.
\]
Consequently,
\[
\max_{(a,b)\in \overline D} J(\gamma(a,b)) < c_R,
\]
hence that $\gamma(\overline D) \cap \mathcal M_R = \emptyset$.  On the other hand, define $\Psi: \overline D \to \mathbb{R}^2$ by
\[
\Psi(a,b) = (a^{-1}J'(\gamma(a,b))(\gamma(a,b)^+), b^{-1}J'(\gamma(a,b))(\gamma(a,b)^-)).
\]
Since $B\cap\{aw_R^+ + bw_R^-\, :\,  (a,b) \in \partial D\} = \emptyset$ and  $\eta(\tau, u) = u$, for every $\tau \geq 0$ provided that $u \not\in B$, we get
\[
\Psi(a,b) = (J'(aw_R^+ + bw_R^-)w_R^+, J'(aw_R^+ + bw_R^-)w_R^-) = \Phi(a,b),\ \forall\,  (a,b) \in \partial D.
\]
By the Brouwer's topological degree,
\[
{\rm{d}}(\Psi, D, (0,0)) = {\rm{d}}(\Phi, D, (0,0)) = {\rm{sgn}}({\rm{det}}(\Phi')(1,1)) = 1,
\]
which implies that there exists $({a_0},{b_0}) \in  D$ such that $J'(\gamma(a_0,b_0))(\gamma(a_0,b_0)^\pm) = 0$. Hence, $\gamma(a_0,b_0) \in \mathcal M_R$. This contradicts the fact that $\gamma(\overline D) \cap \mathcal M_R = \emptyset$.  \qed


\begin{thebibliography}{99}





\bibitem{JMAA2010} C.O. Alves, M.A. Souto, S.H.M. Soares, Schr\"odinger-Poisson equations without Ambrosetti-Rabinowitz
condition, J. Math. Anal. Appl. 377 (2011) 584-592.



\bibitem{nodal1} C.O. Alves, M.A. Souto, Existence of least energy nodal solution for a  Schr\"odinger-Poisson system in bounded domains,
Z. Angew. Math. Phys.  (2013).

\bibitem{AM} {A. Ambrosetti} and {R. Ruiz, }  {Multiple bound states for the {S}chr\"odinger-{P}oisson
              problem}, {Commun. Contemp. Math.} {10} {(2008)}
  {391--404}.


   \bibitem{ap}
    {A. Azzollini} and {A. Pomponio, }
   {Ground state solutions for the nonlinear
              {S}chr\"odinger-{M}axwell equations},
{J. Math. Anal. Appl.}
     {345}
{(2008)}
  {90--108}.


\bibitem{BWW}{T. Bartsch, T. Weth} and { M. Willem,} {Partial symmetry of least energy nodal solution to some variational problems}, {Journal D'Analyse Math\'ematique}{ 1 (2005) 1-18}


\bibitem{BW} T. Bartsch and T. Weth, Three nodal solutions of singularly perturbed elliptic equations
on domains without topology, Ann. Inst. H. Poincar\'e Anal. Non
Lin\'eaire 22 (2005) 259-281.











    \bibitem{bf}
    {V. Benci } and {D. Fortunato, }
 {An eigenvalue problem for the Schr\"odinger-Maxwell equations},
{Top. Meth. Nonlinear Anal.} {11} {(1998)}
     {283--293}.


\bibitem{bl} {H. Berestycki} and {P.L. Lions,}
   {Nonlinear scalar field equations, I - existence of a ground state},
    {Arch. Rat. Mech. Analysis},
   {82},
     {(1983)}
   {313--346}.


\bibitem{BM} {O. Bokanowski} and {N.J. Mauser,}{ Local approximation of the Hartree-Fock exchange potential: a deformation approach,
 ${\text M}^3$AS 9 (1999) 941-961.}








\bibitem{Cerami}
 {G. Cerami} and {G. Vaira,}
 {Positive solutions for some non-autonomous
              {S}chr\"odinger-{P}oisson systems},
{J. Differential Equations}
     {248} {(2010)}
  {521--543}.


\bibitem{C}
 {G.M. Coclite,}
 {A multiplicity result for the nonlinear
              {S}chr\"odinger-{M}axwell equations},
{Commun. Appl. Anal.}
     {7}
{(2003)}
  {417--423}.





\bibitem{Daprile1}
 {T. D'Aprile} and {D. Mugnai},
 {Solitary waves for nonlinear {K}lein-{G}ordon-{M}axwell and
              {S}chr\"odinger-{M}axwell equations},
{Proc. Roy. Soc. Edinburgh Sect. A} {134} {(2004)}
  {893--906}.





\bibitem{Daprile2}
 {T. D'Aprile} and {D. Mugnai},
 {Non-existence results for the coupled  {K}lein-{G}ordon-{M}axwell  equations},
{Adv. Nonlinear Stud.} {4} {(2004)}
  {307--322}.





\bibitem{Davenia}
 {P. d'Avenia,}
 {Non-radially symmetric solutions of nonlinear {S}chr\"odinger
              equation coupled with {M}axwell equations},
{Adv. Nonlinear Stud.}   {2} {(2002)} {177--192}.




     \bibitem{G}
    {G. Siciliano,}
 {Multiple positive solutions for a  Schr\"odinger-Poisson-Slater system},
{J. Math. Analysis and Appl.} {365} {(2010)}
     {288--299}.

\bibitem{K1}
 {H. Kikuchi},
 {On the existence of a solution for elliptic system related to
              the {M}axwell-{S}chr\"odinger equations},
{Nonlinear Anal.}     {67} {(2007)} {1445--1456}.

\bibitem{Ianni2}
 {I. Ianni}, {Sign-Changing radial solutions for the Schr\"odinger-Poisson-Slater problem}, Topol. Methods Nonlinear Anal.  41  (2013) 365–385.

\bibitem{Ianni}
 {I. Ianni} and {G. Vaira},
 {On concentration of positive bound states for the
              {S}chr\"odinger-{P}oisson problem with potentials},
{Adv. Nonlinear Stud.}     {8} {(2008)}  {573-595}.





\bibitem{Kim}
 {S. Kim } and {J. Seok},
 {On nodal solutions of the Nonlinear Schr\"odinger-Poisson equations},
{Comm. Cont. Math.}     {14} {(2012)}  {12450041-12450057}.

\bibitem{lww}
     {Liu, J.-Q}, {Wang, Y.} and {Wang, Z.-Q.}
     {Solutions for quasilinear Schr\"odinger equations via Nehari Method},
    {Comm. in Partial Diff. Equations}
   {29}
     {(2004)}
   {879-901}.






\bibitem{Mauser} {N.J. Mauser,}  {The Schr\"odinger-Poisson-$X_\alpha$ equation,} {Applied Math. Letters 14 (2001)}
759-763.


\bibitem{Miranda} {C. Miranda}, {Un'osservazione su un teorema di Brouwer,}{ Bol. Un. Mat. Ital. 3 (1940)} 5-7.

\bibitem{nehari} Z. Nehari, {Characteristic values associated with a class of nonlinear second order differential equations}, Acta Math. 105 (1961), 141--175.

\bibitem{LG} {L. Pisani} and {G. Siciliano,} { Note on a Schr\"odinger-Poisson system in a bounded domain,
Appl. Math. Lett.} {21} {(2008)}  {521-528} .







     \bibitem{Ruiz}
    {D. Ruiz, }
 {The Schr\"odinger-Poisson equation under the effect of a nonlinear local term },
{J. Funct. Analysis} {237} {(2006)}  {655--674}.



\bibitem{Ruiz-Sic}
 {D. Ruiz} and {G. Siciliano,}
 {A note on the Schr\"odinger-Poisson-{S}later equation on
              bounded domains}, {Adv. Nonlinear Stud.}  {8}
{(2008)}  {179--190}.





\bibitem{Sanchez}
{O. S\'anchez} and {J. Soler}, {Long-time dynamics of the Schr\"odinger-Poisson-Slater system,
J. Statistical Physics 114 (2004) 179-204.}


\bibitem{terracininodea} S. Terracini and G. Verzini, {Solutions of prescribed number of zeroes to a class of superlinear ODE's systems}, NoDEA 8 (2001) 323--341.

































\end{thebibliography}
\end{document}